\documentclass[12pt]{article}

\usepackage{amsmath}
\usepackage{amsthm}
\usepackage{graphicx}
\usepackage{amsfonts}
\usepackage{amssymb}
\usepackage{psfrag}
\usepackage{wrapfig}
\usepackage{url}
\usepackage{amscd}
\usepackage{multicol}
\usepackage{comment}

\newcommand{\re}{\mathbb{R}}

\newcommand{\C}{\mathbb{C}}
\newcommand{\Z}{\mathbb{Z}}

\newcommand{\Q}{\mathbb{Q}}
\renewcommand{\P}{\mathbb{P}}

\newcommand{\dt}{\delta}

\newcommand{\bdes}{\begin{description}}
\newcommand{\edes}{\end{description}}

\newcommand{\bal}{\begin{align}}
\newcommand{\eal}{\end{align}}

\newcommand{\bnum}{\begin{enumerate}}
\newcommand{\enum}{\end{enumerate}}

\newcommand{\bit}{\begin{itemize}}
\newcommand{\eit}{\end{itemize}}

\newcommand{\bea}{\begin{eqnarray*}}
\newcommand{\eea}{\end{eqnarray*}}

\newcommand{\be}{\begin{equation}}
\newcommand{\ee}{\end{equation}}

\newcommand{\baray}{\begin{array}}
\newcommand{\earay}{\end{array}}

\newcommand{\bca}{\begin{cases}}
\newcommand{\eca}{\end{cases}}

\newcommand{\bcen}{\begin{center}}
\newcommand{\ecen}{\end{center}}

\newcommand{\bbm}{\begin{bmatrix}}
\newcommand{\ebm}{\end{bmatrix}}

\newcommand{\bmx}{\begin{matrix}}
\newcommand{\emx}{\end{matrix}}

\newcommand{\bpm}{\begin{pmatrix}}
\newcommand{\epm}{\end{pmatrix}}

\newcommand{\btab}{\begin{tabular}}
\newcommand{\etab}{\end{tabular}}

\theoremstyle{plain}
\newtheorem{thm}{Theorem}[section]
\newtheorem{theorem}[thm]{Theorem}

\newtheorem{lemma}[thm]{Lemma}

\newtheorem{corollary}[thm]{Corollary}

\theoremstyle{definition}

\usepackage[top=1.5in,bottom=1in,left=1.2in,right=1.2in]{geometry}

\begin{document}

\title{Semidefinite Representation of the $k$-Ellipse}

\author{ Jiawang Nie\thanks{Institute for Mathematics and its
Applications, University of Minnesota} \and Pablo
A. Parrilo\thanks{Laboratory for Information and Decision Systems,
Massachusetts Institute of Technology} \and Bernd
Sturmfels\thanks{Department of Mathematics, University of California at
Berkeley}}

\maketitle

\begin{abstract}
\noindent The $k$-ellipse is the plane algebraic curve consisting of
all points whose sum of distances from $k$ given points is a fixed
number.  The polynomial equation defining the $k$-ellipse has degree $2^k$ if $k$
is odd and degree $2^k{-}\binom{k}{k/2}$ if $k$ is even.  We express
this polynomial equation as the determinant of a symmetric matrix of linear
polynomials.  Our representation extends to weighted $k$-ellipses and $k$-ellipsoids in
arbitrary dimensions, and it leads to new geometric applications of
semidefinite programming.
\end{abstract}

\section{Introduction}

The {\em circle} is the plane curve consisting
of all points $(x,y)$ whose distance from a
given point $(u_1,v_1)$ is a fixed number $d$.
It is the zero set of the quadratic polynomial
\be \label{det:circle} p_1(x,y) \quad = \quad {\rm det}  \bbm
                    d+ x  - u_1  &    y - v_1    \\
                        y - v_1    &  d -x  + u_1 \\
\ebm \! .
\ee
The {\em ellipse} is the  plane curve consisting
of all points $(x,y)$ whose sum of distances from two
given points $(u_1,v_1)$ and $(u_2,v_2)$ is a fixed number $d$.
It is the zero set of
\be \label{det:ellipse} \! p_2(x,y) \,\,= \,\, {\rm det}  \! \bbm
d+  2 x - u_1 - u_2 \! &  y - v_1 & y - v_2 & 0 \\
               y - v_1 &    d + u_1 - u_2 &   0 &   y - v_2 \\
               y - v_2 &    0 &   d - u_1 + u_2 &   y - v_1 \\
               0 & y - v_2 & y - v_1 & \! d -2 x + u_1 + u_2 \\
\ebm \! . \ee
In this paper we generalize these determinantal formulas
for the circle and the ellipse.
  Fix a positive real number $d$ and fix $k$ distinct
points $\, (u_1,v_1), (u_2,v_2), \ldots, (u_k,v_k)\,$ in $\re^2$.  The
{\em $k$-ellipse} with {\em foci} $(u_i,v_i)$ and {\em radius} $d$ is
the following curve in the plane:
\be \label{kellipse} \left\{(x,y)
\in \re^2 \,:\,\sum_{i=1}^k \sqrt{(x-u_i)^2+(y-v_i)^2} \, = \,
d\right\}.  \ee

The $k$-ellipse is the boundary of a convex set $\mathcal{E}_k$ in the
plane, namely, the set of points whose sum of distances to the $k$
given points is at most $d$. These convex sets are of interest in
computational geometry \cite{CINDY} and in optimization, e.g.~for the
Fermat-Weber facility location problem \cite{Bajaj,CT,
Kulshrestha,Sekino,Weiszfeld}.  In the classical literature
(e.g.~\cite{Sturm}), $k$-ellipses are known as {\em Tschirnhaus'sche
Eikurven} \cite{Nagy}. Indeed, they look like ``egg curves'' and they
were introduced by Tschirnhaus in 1686.

We are interested in the irreducible polynomial $p_k(x,y)$ that
vanishes on the $k$-ellipse.  This is the unique (up to sign)
polynomial with integer coefficients in the unknowns $x$ and $y$ and
the parameters $d,\,u_1,v_1,\,\ldots,u_k,v_k$.  By the {\em degree of
the $k$-ellipse} we mean the degree of $p_k(x,y)$ in $x$ and $y$.  To
compute it, we must eliminate the square roots in the representation
(\ref{kellipse}).  Our solution to this problem is as follows:

\begin{theorem} \label{thm1}
The $k$-ellipse has degree $2^k$ if $k$ is odd and degree
$2^k{-}\binom{k}{k/2}$ if $k$ is even.  Its defining polynomial has a
determinantal representation \be \label{detrep} p_k(x,y) \quad = \quad
{\rm det}\bigl(\,x \cdot A_k + y \cdot B_k + C_k \,\bigr) \ee where
$A_k, B_k,C_k$ are symmetric $2^k \times 2^k$ matrices.  The entries
of $A_k$ and $B_k$ are integer numbers, and the entries of $C_k$ are
linear forms in the parameters $d,u_1,v_1,\ldots,u_k,v_k$.
\end{theorem}

\begin{figure}
\centering
\includegraphics[width=4.7cm]{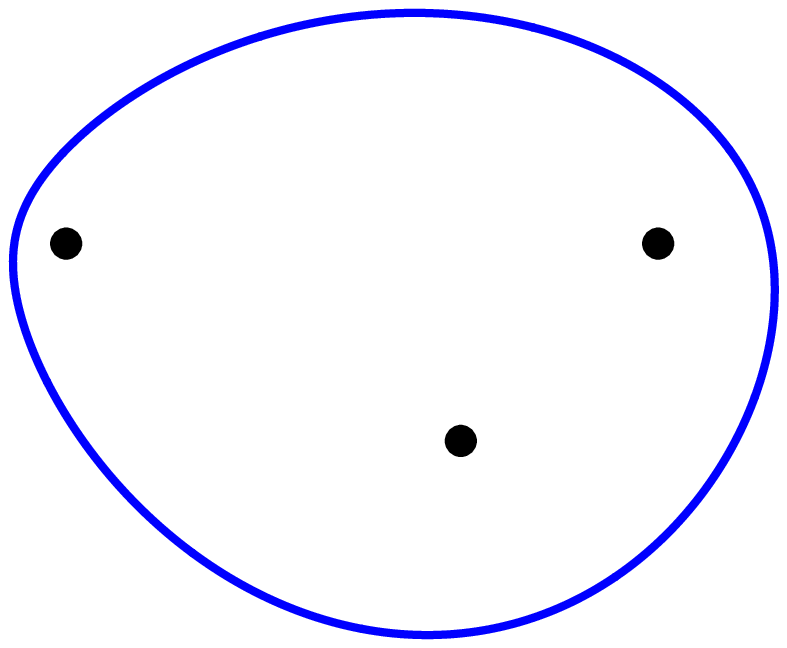}\quad \,\,\,
\includegraphics[width=4.6cm]{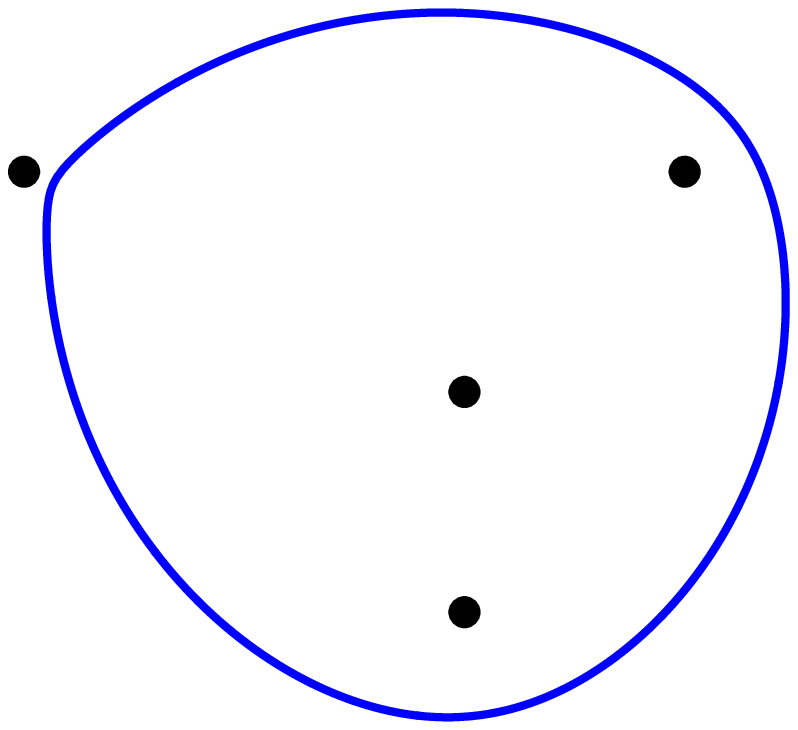}\quad \,\,\,
\includegraphics[width=4.7cm]{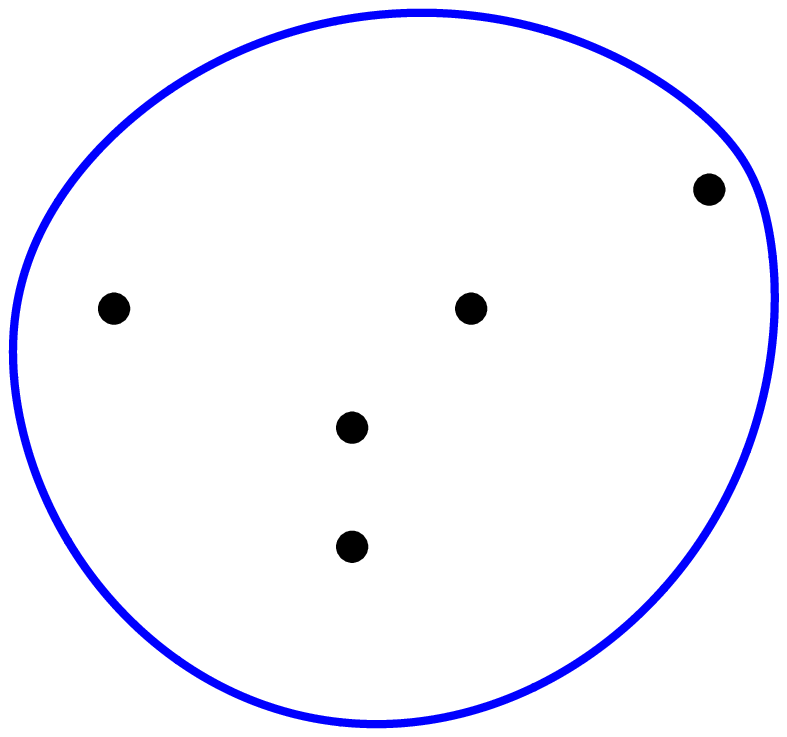}
\caption{A $3$-ellipse, a $4$-ellipse, and a $5$-ellipse, each with its foci.}
\label{fig:three}
\end{figure}

For the circle $(k=1)$
and the ellipse $(k=2)$, the representation (\ref{detrep})
is given by the formulas (\ref{det:circle}) and (\ref{det:ellipse}).
The polynomial $p_3(x,y)$ for the $3$-ellipse is the determinant of
{\tiny
$$ \bbm
d {+} 3 x  {-} u_1 {-} u_2 {-} u_3 \!\!\!\!\!\!\! &  y {-} v_1 & y {-} v_2 & 0 & y {-} v_3 & 0 & 0 & 0 \\
  y {-} v_1 & \!\!\!\!\!\!\! d {+} x  {+} u_1 {-} u_2 {-} u_3 \!\!\!\!\!\!\! &
  0 &  y {-} v_2 & 0 & y {-} v_3 & 0 & 0 \\
  y {-} v_2 &  0 & \!\!\!\!\!\!\! d {+} x  {-} u_1 {+} u_2 {-} u_3 \!\!\!\!\!\!\! &
  y {-} v_1 & 0 & 0 & y {-} v_3 & 0 \\
  0 &  y {-} v_2 &  y {-} v_1 & \!\!\!\!\!\!\! d {-}x  {+} u_1 {+} u_2 {-} u_3
\!\!\!\!\!\!\! & 0 & 0 & 0 & y {-} v_3 \\
  y {-} v_3 &  0 &  0 &  0 & \!\!\!\!\!\!\! d {+} x {-} u_1 {-} u_2 {+} u_3
\!\!\!\!\!\!\! & y {-} v_1 & y {-} v_2 & 0 \\
  0 &  y {-} v_3 &  0 & 0 & y {-} v_1 & \!\!\!\!\!\!\!
d {-}x  {+} u_1 {-} u_2 {+} u_3 & 0 & y {-} v_2 \!\!\!\!\!\!\! \\
  0 &  0 &  y {-} v_3 & 0 & y {-} v_2 & 0 &
\!\!\!\!\!\!\! d {-}x  {-} u_1 {+} u_2 {+} u_3 & y {-} v_1 \!\!\!\!\!\!\! \\
  0 & 0 & 0 & y {-} v_3 & 0 & y {-} v_2 & y {-} v_1 &
\!\!\!\!\!\!\! d {-}3 x {+} u_1 {+} u_2 {+} u_3  \\
\ebm $$
}
The full expansion of this $8 \times 8$-determinant has $2,355$ terms.
Next, the $4$-ellipse is a curve of degree ten which is
represented by a symmetric $16 \times 16$-matrix, etc....

\begin{figure}
\centering
\includegraphics[height=11cm]{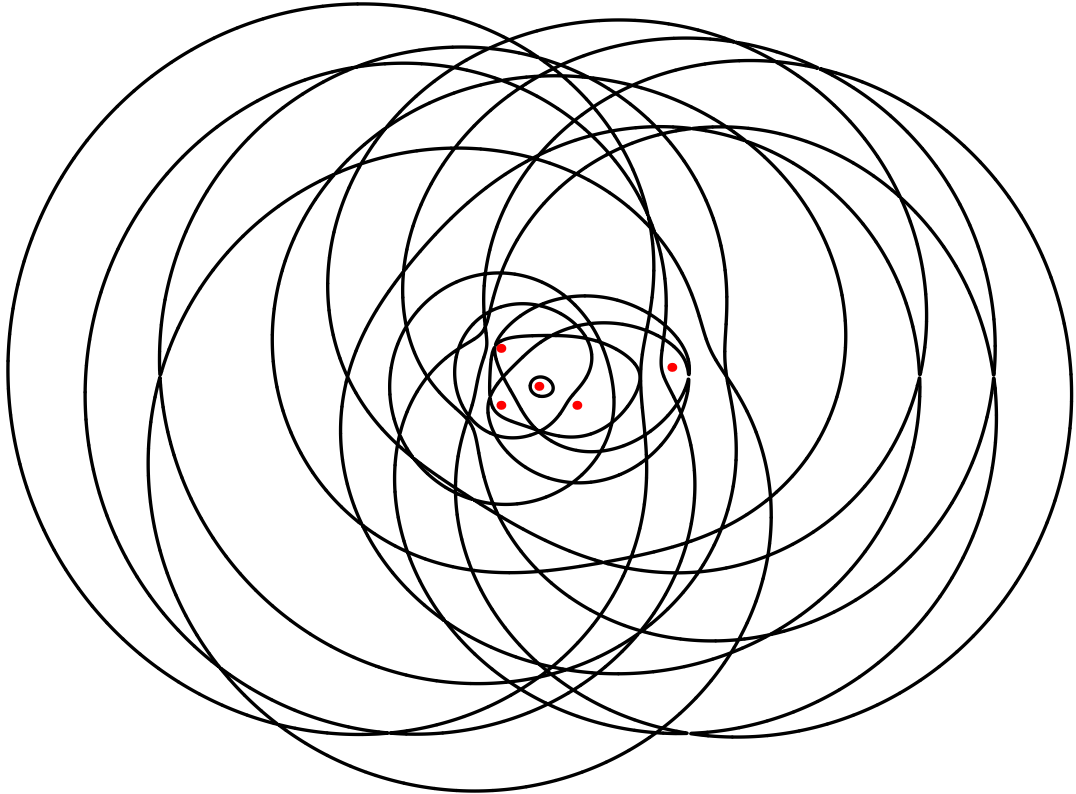}
\vskip -0.4cm
\caption{The Zariski closure of the 5-ellipse is an algebraic
curve of degree 32.}
\label{fig:components5ellipse}
\end{figure}

This paper is organized as follows.  The proof of Theorem \ref{thm1}
will be given in Section 2.  Section 3 is devoted to  geometric aspects
and connections to semidefinite programming.
  While the $k$-ellipse itself is a convex curve,
its Zariski closure $\,\{\, \,p_k(x,y) = 0 \,\}\,$ has many extra
branches outside the convex set $\,\mathcal{E}_k$. They are arranged
in a beautiful pattern known as a {\em Helton-Vinnikov curve} \cite{HV}.
This pattern is shown in Figure
\ref{fig:components5ellipse} for $k=5$ points.
In Section 4 we generalize our results
to higher dimensions and to the weighted case, and we discuss the
computation of the {\em Fermat-Weber point} of the given points
$(u_i,v_i)$.  A list of open problems and future directions is presented in Section~5.

\section{Derivation of the Matrix Representation}

We begin with a discussion of the degree of the $k$-ellipse.

\begin{lemma} \label{lem21}
The defining polynomial of the $k$-ellipse has degree at most $2^k$ in
the variables $(x,y)$ and it is monic of degree $2^k$ in the radius
parameter~$d$.
\end{lemma}

\begin{proof}
We claim that the defining  polynomial
of the $k$-ellipse can be written as follows:
\begin{equation}
\label{galoisproduct}
 p_k(x,y) \quad = \quad
 \prod_{\sigma \in \{0,1\}^k }
\left( \,d \,-\, \sum_{i=1}^k (-1)^{\sigma_i}
\sqrt{(x-u_i)^2+(y-v_i)^2} \,\,\right).
\end{equation}
Obviously, the right hand side vanishes on the $k$-ellipse.
The following Galois theory argument shows that this expression
is a  polynomial and that it is irreducible.
Consider the polynomial ring $\, R =  \Q[x,y,d,u_1,v_1,\ldots,u_k,v_k]$.
The field of fractions of $R$ is the field
$\,K =  \Q(x,y,d,u_1,v_1,\ldots,u_k,v_k)\,$
of rational functions in all unknowns.
Adjoining the square roots in (\ref{kellipse}) to $K$ gives
 an algebraic field extension $L$ of degree $2^k$ over $K$.
 The Galois group of the extension $L/K$ is $(\Z/2\Z)^k$, and
 the product in (\ref{galoisproduct}) is over the orbit of the element
$\,d- \sum_{i=1}^k \sqrt{(x-u_i)^2+(y-v_i)^2}\,$ of $L$
under the action of the Galois group. Thus this product
in (\ref{galoisproduct}) lies in the ground field $K$.
Moreover, each factor in the product is integral over $R$,
and therefore the product lies in the polynomial ring $R$.
To see that this polynomial is irreducible, it suffices to observe
that no proper subproduct of the right hand side in (\ref{galoisproduct}) lies in
the ground field $K$.
\end{proof}

The statement {\em degree at most $2^k$} is the crux in Lemma \ref{lem21}.
Indeed, the degree in $(x,y)$ and can be strictly smaller than $2^k$
as the case of the classical ellipse $(k=2)$ demonstrates.  When
evaluating the product (\ref{galoisproduct}) some unexpected
cancellations may occur.  This phenomenon happens for all even $k$,
as we shall see later in this section.

We now turn to the matrix representation promised by
Theorem~\ref{thm1}. We recall the following standard definition
from matrix theory (e.g., \cite{HJ2}).
Let $A$ be a real $m \times m$-matrix and
$B$ a real $n \times n$-matrix.  The {\em tensor
 sum} of $A$ and $B$ is the $mn \times mn$ matrix $A
\oplus B := A \otimes I_n + I_m \otimes B $. The tensor sum of square matrices
is an associative operation which is not commutative.
For instance, for three matrices $A, B,C$ we have
$$ A \oplus B \oplus C \quad = \quad
A \otimes I \otimes I  \, + \,
I \otimes B \otimes I \, + \,
I \otimes I \otimes C . $$
Here $\otimes$ denotes the {\em tensor product},
which is also associative but not commutative.
Tensor products and tensor sums of matrices are also known as
\emph{Kronecker products} and
\emph{Kronecker sums} \cite{Bellman,HJ2}.
Tensor sums of symmetric matrices can be diagonalized
by treating the summands separately:

\begin{lemma} \label{lem:diagonalization}
Let $M_1,\ldots,M_k$ be symmetric matrices,
$U_1,\ldots,U_k$ orthogonal matrices,
and $\Lambda_1,\ldots,\Lambda_k$ diagonal matrices
such that $\,M_i = U_i \cdot \Lambda_i \cdot U_i^T\,$ for $ i=1,\ldots,k$. Then
\[
(U_1 \otimes \cdots \otimes U_k)^T \cdot
(M_1 \oplus \cdots \oplus M_k) \cdot
(U_1 \otimes \cdots \otimes U_k) \,\,=\,\,
\Lambda_1 \oplus \cdots \oplus \Lambda_k.
\]
In particular, the eigenvalues of the tensor sum $M_1 \oplus M_2 \oplus \cdots \oplus M_k$
are the sums $\lambda_1 + \lambda_2 + \cdots + \lambda_k$ where
$\lambda_1$ is any eigenvalue of $M_1$,
$\lambda_2$ is any eigenvalue of $M_2$, etc.
\end{lemma}

The proof of this lemma is an exercise in (multi)-linear algebra.
We are now prepared to state our formula for the
 explicit determinantal representation of the $k$-ellipse.

\begin{thm} \label{lem22}
Given points $(u_1,v_1),\ldots,(u_k,v_k)$ in $\re^2$, we define the
 $2^k \times 2^k$ matrix
\be
L_{k}(x,y) \,\,\,:= \,\,\, d \cdot I_{2^k} \,+\,
\bbm x - u_1 & y -v_1 \\ y -v_1 & - x + u_1 \ebm \oplus \,\cdots\, \oplus
\bbm x - u_k & y -v_k \\ y -v_k & - x + u_k \ebm
\label{eqn:LABC}
\ee
which is affine in $x,y$ and $d$. Then the $k$-ellipse has the
determinantal representation
\begin{equation}
p_k(x,y) \,\,= \,\, \det L_k(x,y).
\label{eq:tensorrep}
\end{equation}
The convex region bounded by the $k$-ellipse
is defined by the following matrix inequality:
\begin{equation}
\mathcal{E}_k \,\,= \,\,
\left\{ \,(x,y) \in \re^2 \,:\,
L_k(x,y) \,\,\succeq\, \,0 \,\right\}.
\end{equation}
\end{thm}

\begin{proof}
Consider the $2 \times 2$ matrix that appears
as a tensor summand in~(\ref{eqn:LABC}):
\[
\bbm x - u_i & y -v_i \\ y -v_i & - x + u_i \ebm.
\]
A computation shows that this matrix is
orthogonally similar to
\[
\bbm \sqrt{(x-u_i)^2+(y-v_i)^2} & 0 \\
0 & -\sqrt{(x-u_i)^2+(y-v_i)^2} \ebm.
\]
These computations take place in the field $L$
which was considered in the proof of
Lemma \ref{lem21} above.
Lemma~\ref{lem:diagonalization} is valid over any field, and it
implies that the matrix
$L_k(x,y)$ is orthogonally similar to a $2^k \times 2^k$ diagonal
matrix with diagonal entries
\begin{equation}
d \,+\, \sum_{i=1}^k (-1)^{\sigma_i} \sqrt{(x-u_i)^2+(y-v_i)^2}, \qquad \sigma_i \in \{0,1\}.
\label{eq:eigenvals}
\end{equation}
The desired identity~(\ref{eq:tensorrep}) now follows directly
from~(\ref{galoisproduct}) and the fact that the determinant of a
matrix is the product of its eigenvalues.
 For the characterization of the convex set
$\mathcal{E}_k$, notice that positive semidefiniteness of $L_k(x,y)$
is equivalent to nonnegativity of all the
eigenvalues~(\ref{eq:eigenvals}). It suffices to consider
the smallest eigenvalue, which equals
\[
d \, - \,\sum_{i=1}^k \sqrt{(x-u_i)^2+(y-v_i)^2} .
\]
Indeed, this quantity is nonnegative if and only if
the point $(x,y)$ lies in $ \mathcal{E}_k$.
\end{proof}

\begin{proof}[Proof of Theorem \ref{thm1}]
The assertions in the second and third sentence have just been
proved in Theorem \ref{lem22}.  What remains to be shown is the first
assertion concerning the degree of $\,p_k(x,y)\,$ as a polynomial in
$(x,y)$.  To this end, we consider the univariate polynomial
 $\,g(t)\,:=\,p_k(t \cos \theta, t \sin \theta)$ where
$\theta$ is a generic angle.  We must prove that
\[
{\rm deg}_t \bigl( g(t) \bigr) \,\, = \,\,\bca
\,\, 2^k & \text{ if $k$ is odd, } \\
\,\, 2^k-\binom{k}{k/2} & \text{ if $k$ is even. }
\eca
\]
The polynomial $\,g(t)\,$ is the determinant of the symmetric $2^k \times 2^k$-matrix
\begin{equation}
L_k(t \cos \theta, t \sin \theta) \,\,\, = \,\,\, t \cdot \left(
\bbm \cos \theta & \!\phantom{-} \sin \theta \\ \sin \theta & \! - \cos \theta \ebm \oplus \cdots \oplus
\bbm \cos \theta & \!\phantom{-} \sin \theta \\ \sin \theta & \! - \cos \theta \ebm \right) \,+\, C_k.
\label{eq:lkt}
\end{equation}
The matrix $C_k$ does not depend on $t$.
We now define the $2^k \times 2^k$ orthogonal matrix
\[
\quad U \,\, :=\,\, \underbrace{V \otimes \cdots \otimes V}_{k \; \mathrm{times}} \qquad
\hbox{where} \quad
V \,:=\, \bbm
\cos (\theta/2) & \! -\sin (\theta / 2) \\
\sin  (\theta/2) & \! \phantom{-} \cos (\theta/2) \ebm,
\]
and we note the matrix identity
\[
V^T \cdot
\bbm \cos \theta & \phantom{-} \sin \theta \\
 \sin \theta & - \cos \theta \ebm \cdot
V \quad = \quad \bbm 1 & 0 \\ 0 & -1 \ebm. \]
Pre- and post-multiplying~(\ref{eq:lkt}) by $U^T$ and $U$, we find that
$$ U^T \cdot L_k(t \cos \theta, t \sin \theta) \cdot U
\quad  = \quad  t \cdot
\underbrace{ \left(
\bbm 1 & 0 \\ 0 & -1\ebm \oplus \cdots \oplus
\bbm 1 & 0 \\ 0 & -1 \ebm \right) }_{E_k}  \,+ \,\, U^T \cdot C_k \cdot U .$$
The determinant of this matrix is our univariate polynomial $g(t)$.
The matrix $E_k$ is a diagonal matrix of dimension $2^k \times 2^k$.
Its diagonal entries are obtained by summing $k$
copies of $-1$ or $+1$ in all $2^k$ possible ways. None of these sums
are zero when $k$ is odd, and precisely $\binom{k}{k/2}$ of these sums
are zero when $k$ is even.  This shows that the rank of $E_k$ is $2^k$
when $k$ is odd, and it is $2^k - \binom{k}{k/2}$ when $k$ is even.
We conclude that the univariate polynomial $\,g(t) \, = \,
{\rm det} \bigl(\, t \cdot E_k + U^T C_k U \bigr)\,$ has the
desired degree.
\end{proof}

\begin{figure}
\centering
\includegraphics[height=9.3cm]{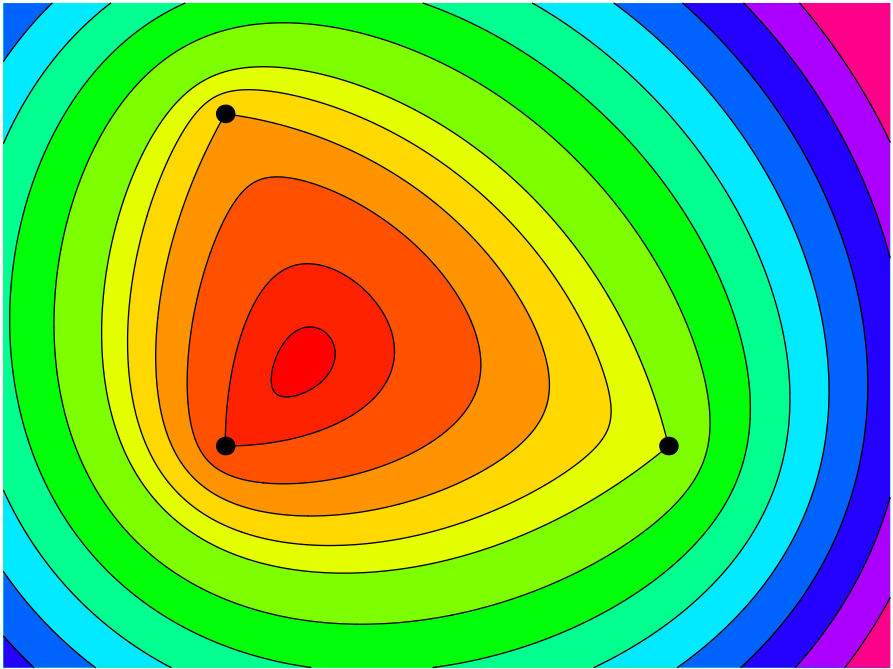}
\caption{A pencil of 3-ellipses with fixed foci (the three black dots)
and different radii.  These convex curves are always smooth
unless they contain one of the foci.}
\label{fig:concentric}
\end{figure}

\section{More Pictures and Some Semidefinite Aspects}
\label{sec:semidefinite}

In this section we examine the geometry of the $k$-ellipse, we
look at some pictures, and we discuss aspects
relevant to the theory of semidefinite programming. In
Figure~\ref{fig:three} several $k$-ellipses are shown, for $k=3,4,5$.
One immediately observes that, in contrast to the classical
circle and ellipse, a $k$-ellipse does not necessarily contain the
foci in its interior.
The interior $\mathcal{E}_k$ of the $k$-ellipse is a sublevel set of the
convex function
\be \label{convexfct} (x,y) \,\,\,\, \mapsto \,\,\,\, \sum_{i=1}^k \sqrt{(x-u_i)^2+(y-v_i)^2}.
\ee
This function is strictly convex
in any region excluding the foci $\{(u_i,v_i)\}_{i=1}^k$,
provided the foci are not collinear \cite{Sekino}.
This explains why the $k$-ellipse is a convex curve. In order for
$\mathcal{E}_k$ to be nonempty, it is necessary
and sufficient that the radius $d$ be
greater than or equal to the global minimum $d_\star$
of the convex function (\ref{convexfct}).

The point $(x_\star,y_\star)$ at which the global minimum $d_\star$ is
achieved is called the \emph{Fermat-Weber point} of the foci.  This
point minimizes the sum of the distances to the $k$ given points
$(u_i,v_i)$, and it is of importance in the facility location
problem. See \cite{Bajaj,CT, Kulshrestha}, and \cite{Sturm} for a
historical reference.  For a given set of foci, we can vary the radius
$d$, and this results in a pencil of confocal $k$-ellipses, as in
Figure~\ref{fig:concentric}. The sum of distances function
(\ref{convexfct}) is differentiable everywhere except at the
$(u_i,v_i)$, where the square root function has a singularity. As a
consequence, the $k$-ellipse is a smooth convex curve, except when
that curve contains one of the foci.

An algebraic geometer would argue that there is more to the
$k$-ellipse than meets the eye in Figures~\ref{fig:three} and
\ref{fig:concentric}. We define the {\em algebraic $k$-ellipse} to be
the Zariski closure of the $k$-ellipse, or, equivalently, the zero set
of the polynomial $p_k(x,y)$. The algebraic $k$-ellipse is an
algebraic curve, and it can be considered in either the real plane
$\re^2$, in the complex plane $\C^2$, or (even better) in
the complex projective plane $\P^2_\C$.

Figure~\ref{fig:components5ellipse} shows an algebraic $5$-ellipse.
In that picture, the actual $5$-ellipse is the tiny convex curve
in the center. It surrounds only
one of the five foci.

\begin{figure}
\centering
\includegraphics[height=9cm]{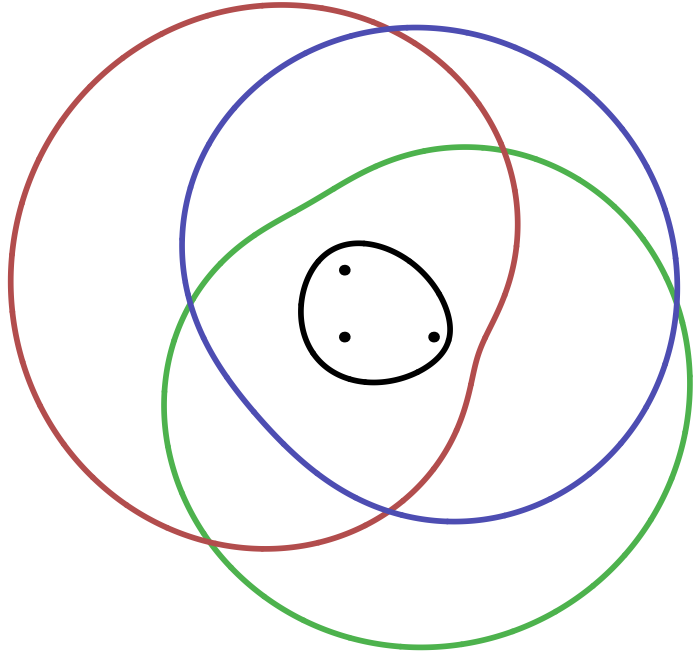}
\vskip -0.4cm
\caption{The Zariski closure of the 3-ellipse is an algebraic
curve of degree eight.}
\label{fig:components3ellipse}
\end{figure}

For a less dizzying illustration see
Figure~\ref{fig:components3ellipse}. That picture shows an algebraic
$3$-ellipse.  The curve has degree eight, and it is given algebraically by the
$8 \times 8$-determinant displayed in the Introduction.
We see that the set of real points on the algebraic $3$-ellipse consists of
four ovals, corresponding to the equations
$$ \sqrt{(x-u_1)^2+(y-v_1)^2}
\pm \sqrt{(x-u_2)^2+(y-v_2)^2}
\pm \sqrt{(x-u_3)^2+(y-v_3)^2} \,\,\, = \,\,\,  d .$$
Thus Figure~\ref{fig:components3ellipse} visualizes
the Galois theory argument  in the proof of Lemma \ref{lem21}.

\begin{figure}
\begin{center}
\includegraphics[height=9.7cm]{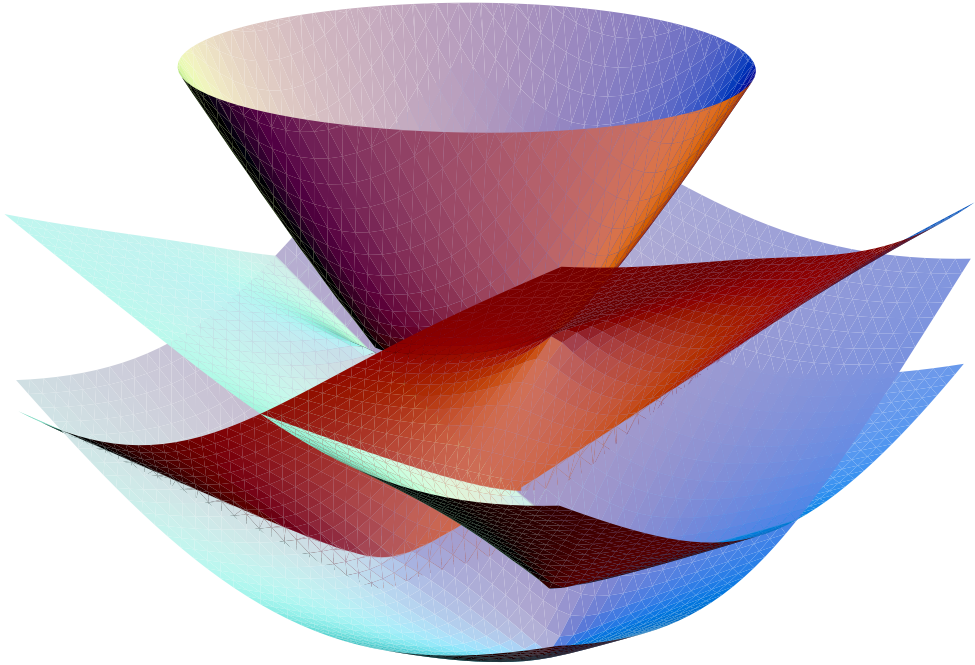}
\end{center}
\vskip -0.9cm
\caption{The irreducible surface $\{\,p_3(x,y,d)=0\}$. Taking
horizontal slices gives the pencil of algebraic $3$-ellipses
for fixed foci and different radii $d$, as shown in Figure \ref{fig:concentric}.}
\label{fig:levelsets}
\end{figure}

If we regard the radius $d$ as an unknown, in addition to the two
unknowns $x$ and $y$, then the determinant in Theorem \ref{thm1}
specifies an irreducible surface $\,\{p_k(x,y,d) = 0\} \,$ in
three-dimensional space.  That surface has degree $2^k$. For an
algebraic geometer, this surface would live in complex projective
$3$-space $\C \P^3$, but we are interested in its points in real
affine $3$-space $\re^3$.  Figure~\ref{fig:levelsets} shows this
surface for $k=3$.  The bowl-shaped convex branch near the top is the
graph of the sum of distances function (\ref{convexfct}), while each
of the other three branches is associated with a different combination
of signs in the product~(\ref{galoisproduct}).  The surface has a
total of $2^k=8$ branches, but only the four in the half-space $d \geq
0$ are shown, as it is symmetric with respect to the plane $\,d=0$.
Note that the Fermat-Weber point $(x_*,y_*,d_*)$ is a highly singular
point of the surface.

 The time has now come for us to explain the adjective
``semidefinite'' in the title of this paper. {\em Semidefinite programming} (SDP)
is a widely used method in convex optimization.
Introductory references include \cite{SIAM,Handbook}.
An algebraic perspective was recently given in \cite{NRS}.
The problem of SDP is
to minimize a linear functional over the solution set of a linear matrix inequality (LMI).
An example of an LMI is
\be \label{ourLMI}
 x \cdot A_k + y \cdot B_k + d \cdot I_{2^k} + \tilde{C}_k \,\,\succeq \,\, 0 .
 \ee
Here $\tilde{C}_k$ is the matrix gotten from $C_k$ by setting $d=0$, so
 that $\,C_k = d \cdot I_{2^k} + \tilde{C}_k$. If $d$ is a fixed positive real number then
the solution set to the LMI  (\ref{ourLMI}) is the convex region
$\mathcal{E}_k$ bounded by the $k$-ellipse. If $d$ is an unknown then
the solution set to the LMI (\ref{ourLMI}) is the epigraph of (\ref{convexfct}),
or, geometrically, the unbounded
$3$-dimensional convex region interior to the bowl-shaped
surface in Figure~\ref{fig:levelsets}. The bottom point
of that convex region is the Fermat-Weber point $(x_*,y_*,d_*)$,
and it can be computed by solving the SDP
\be \label{FWSDP}
 {\rm Minimize} \,\, d \,\, \hbox{subject to} \,\, (\ref{ourLMI}). \ee
Similarly, for fixed radius $d$, the $k$-ellipse is the set of all solutions to
\be
\label{ellipseSDP}
{\rm Minimize} \,\,\,\alpha x + \beta y \,\, \,
\hbox{subject to} \,\, (\ref{ourLMI})
\ee
where $\alpha,\beta$ run over $\re$.
This explains the term {\em semidefinite representation} in our title.

While the Fermat-Weber SDP (\ref{FWSDP}) has only three unknowns, it
has the serious disadvantage that the matrices are
exponentially large (size $2^k$). For computing $(x_*,y_*,d_*)$ in
practice, it is better to introduce slack variables $d_1,d_2,\ldots,d_k$,
and to solve
\be
\label{betterFWSDP}
{\rm Minimize} \,\,\sum_{i=1}^k d_i \,\,\,\hbox{subject to} \,\,\,
\bbm d_i + x-a_i & y-b_i \\ y-b_i & d_i -x+a_i \ebm \, \succeq \,0 \quad (i = 1,\ldots, k).
\ee
This system can be written as single LMI by stacking the $2 \times
2$-matrices to form a block matrix of size $2k \times 2k$. The size of
the resulting LMI is linear in $k$ while the size of the LMI
(\ref{FWSDP}) is exponential in $k$.  If we take the LMI
(\ref{betterFWSDP}) and add the linear constraint $\,d_1 + d_2 +
\cdots + d_k = d$, for some fixed $d > 0$, then this furnishes a
natural and concise \emph{lifted} semidefinite representation of our
$k$-ellipse. Geometrically, the representation expresses
$\mathcal{E}_k$ as the \emph{projection} of a convex set defined by
linear matrix inequalities in a higher-dimensional space.  Theorem
\ref{thm1} solves the algebraic {\em LMI elimination problem}
corresponding to this projection, but at the expense of an exponential
increase in size,
which is due to the exponential growing of degrees of $k$-ellipses.

\smallskip

Our last topic in this section is the relationship of the $k$-ellipse
to the celebrated work of Helton and Vinnikov \cite{HV} on LMI
representations of planar convex sets, which led to the resolution of
the Lax conjecture in~\cite{LPR}.  A semialgebraic set in the plane is
called {\em rigidly convex} if its boundary has the property that
every line passing through its interior intersects the Zariski closure
of the boundary only in real points.  Helton and Vinnikov
\cite[Thm.~2.2]{HV} proved that a plane curve of degree $d$ has an LMI
representation by symmetric $d \times d$ matrices if and only if the
region bounded by this curve is rigidly convex.  In arbitrary
dimensions, rigid convexity holds for every region bounded by a
hypersurface that is given by an LMI representation, but the strong
form of the converse, where the degree of the hypersurface precisely
matches the matrix size of the LMI, is only valid in two dimensions.

It follows from the LMI representation  in
Theorem~\ref{thm1} that the region bounded by a $k$-ellipse is rigidly convex.
Rigid  convexity can be  seen in Figures~\ref{fig:components3ellipse}
and \ref{fig:components5ellipse}.
Every  line  that  passes  through  the interior  of  the  $3$-ellipse
intersects the algebraic $3$-ellipse  in eight real points,
and lines through the $5$-ellipse meet its Zariski closure
in  $32$ real points.  Combining our Theorem \ref{thm1}
with the Helton-Vinnikov Theorem, we conclude:

\begin{corollary}
The $k$-ellipse is rigidly convex. If $k$ is odd, it can be
represented by an LMI of size $2^k$, and if $k$ is even, it can be
represented as by LMI of size $2^k-\binom{k}{k/2}$.
\label{cor:lmis}
\end{corollary}

We have not found yet an \emph{explicit} representation of size
$2^k-\binom{k}{k/2}$ when $k$ is even and $k \geq 4$. For the
classical ellipse $(k=2)$, the determinantal representation
(\ref{det:ellipse}) presented in the Introduction has size $4 \times
4$, while Corollary~\ref{cor:lmis} guarantees the existence of a $2
\times 2$ representation. One such representation of the ellipse with
foci $(u_1,v_1)$ and $(u_2,v_2)$ is: {\small
\[
\bigl( d^2+(u_1{-}u_2)(2x{-}u_1{-}u_2) + (v_1{-}v_2)(2y{-}v_1{-}v_2) \bigr)
\cdot I_2  \,+\,
2d \cdot \bbm x-u_2  & y-v_2 \\ y-v_2 & -x+u_2  \ebm
\,\,\succeq \,\, 0.
\]
} Notice that in this LMI representation of the ellipse, the matrix entries are linear
in $x$ and $y$, as required, but they are
 quadratic in the radius parameter $d$ and the foci $u_i,v_i$.
What is the nicest generalization of this representation
to the $k$-ellipse for $k$ even~?

\section{Generalizations}

The semidefinite representation of the $k$-ellipse we have found
in Theorem \ref{thm1} can
be generalized in several directions. Our first generalization
corresponds to the inclusion of arbitrary nonnegative weights for the
distances, while the second one extends the results from plane curves
to higher dimensions. The resulting geometric shapes are known
{\em Tschirnhaus'sche Eifl\"achen} (or ``egg surfaces'') in the
 classical literature \cite{Nagy}.

\subsection{Weighted $k$-ellipse}

Consider $k$ points $(u_1,v_1),\ldots,(u_k,v_k)$ in the real plane
$\re^2$, a positive radius $d$, and positive \emph{weights}
$w_1,\ldots,w_k$. The {\it weighted $k$-ellipse} is the plane curve
defined as
\[
\left \{(x,y) \in \re^2\; : \,\sum_{i=1}^k w_i \cdot \sqrt{(x-u_i)^2+(y-v_i)^2} \, = \, d \, \right \},
\]
where $w_i$ indicates the relative weight of the distance from $(x,y)$
to the $i$-th focus $(u_i,v_i)$. The \emph{algebraic weighted
$k$-ellipse} is the Zariski closure of this curve.
It is the zero set of an irreducible
polynomial $\,p_k^w(x,y)\,$ that can be constructed as in
equation~(\ref{galoisproduct}). The interior of the weighted
$k$-ellipse is the bounded convex region
\[
\mathcal{E}_k(w) \,\, := \,\, \left\{(x,y) \in \re^2 \, : \, \sum_{i=1}^k  w_i \cdot
\sqrt{(x-u_i)^2+(y-v_i)^2} \, \leq \, d \,\right\}.
\]
The matrix construction from the unweighted
case in (\ref{eqn:LABC}) generalizes as follows:
\be
L^w_{k}(x,y) \,\,\,:= \,\,\, d \cdot I_{2^k} \,+\,
w_1 \cdot
\bbm x - u_1 & y -v_1 \\ y -v_1 & - x + u_1 \ebm \oplus \,\cdots\, \oplus
\, w_k \cdot \bbm x - u_k & y -v_k \\ y -v_k & - x + u_k \ebm
\label{eqn:LABCw}
\ee
Each tensor summand is simply multiplied by the
corresponding weight.
The following representation theorem and degree formula
are a direct generalization of Theorem \ref{lem22}:

\begin{thm}
The algebraic weighted $k$-ellipse
has the semidefinite representation
\[
p_k^w(x,y) \,\,= \,\,\det L^w_{k}(x,y),
\]
and the convex region in its interior satisfies
\[
\mathcal{E}_k(w) \,\,\,= \,\,\,\left\{\, (x,y) \in \re^2 \, :\,
L^w_{k}(x,y)\, \succeq \, 0 \,\right\}.
\]
The degree of the weighted $k$-ellipse is given by
\[
\deg p_k^w(x,y) \,\,= \,\, 2^k - |\mathcal{P}(w)|,
\]
where $\, \mathcal{P}(w) = \{ \dt\in \{-1,1\}^k: \sum_{i=1}^k \dt_i w_i =0 \}$.
\label{thm:weights}
\end{thm}

\begin{proof}
The proof is entirely analogous to that of  Theorem~\ref{lem22}:
\end{proof}

A consequence of the characterization above is the
following cute complexity result.

\begin{corollary}
The decision problem  ``Given a
weighted $k$-ellipse with fixed foci and positive integer
weights, is its algebraic degree smaller than
 $2^k$\,$?$'' is NP-complete.
\end{corollary}

\begin{proof}
Since the number partitioning problem is NP-complete
\cite{GareyJohnson},
this follows from  the degree formula in Theorem~\ref{thm:weights}.

\end{proof}

\subsection{$k$-Ellipsoids}

The definition of a $k$-ellipse in the plane can be naturally extended
to a higher-dimensional space to obtain $k$-ellipsoids.
Consider $k$ fixed points ${\bf u}_1 ,\ldots,{\bf u}_k$ in $\re^n$,
with ${\bf u}_i = (u_{i1}, u_{i2}, \ldots,u_{in})$,
the {\em $k$-ellipsoid} in $\re^{n}$ with these foci is the hypersurface
\be
\label{eiflache1}
 \bigl\{ \,{\bf x} \in  \re^n \,:\,\,
\sum_{i=1}^k \| {\bf u}_i - {\bf x} \| \,= \, d \, \bigr\} \,\,\, = \,\,\,
\left\{ \, {\bf x} \in  \re^n \,:\,\,
\sum_{i=1}^k \sqrt{ \sum_{j=1}^n (u_{ij} - x_j)^2}  \,\,= \,\, d \,\right\}.
\ee
This hypersurface encloses the convex region
\[
\mathcal{E}_k^n \quad = \quad \bigl\{ \, {\bf x} \in  \re^n \, \,:\,\,
\sum_{i=1}^k  \| {\bf u}_i - {\bf x} \| \,\leq \, d \, \bigr\}.
\]
The Zariski closure of the $k$-ellipsoid is
defined by an irreducible polynomial
$\,p_k^n({\bf x}) = p_k^n(x_1,x_2,\ldots,x_n)$.
By the same reasoning as in Section 2, we can prove
the following:

\begin{thm} \label{thm:ellipsoid}
The defining irreducible polynomial $p^k_n({\bf x})$ of the
$k$-ellipsoid is  monic of degree $2^k$ in the parameter~$d$,
it has degree $2^k$ in ${\bf x}$ if $k$ is odd, and it has degree $2^k
- \binom{k}{k/2}$ if $k$ is even.
\end{thm}

We shall represent the polynomial $p^k_n({\bf x})$ as
a factor of the determinant of a symmetric matrix of
affine-linear forms.
To construct this semidefinite representation of the $k$-ellipsoid, we
proceed as follows.  Fix an integer $m \geq 2$.  Let $\,{\bf U}_i({\bf
x})$ be any symmetric $m {\times} m$-matrix of rank $2$ whose entries
are affine-linear forms in ${\bf x}$, and whose two non-zero eigenvalues are
$\,\pm \| {\bf u}_i - {\bf x} \| $.  Forming the tensor sum of these matrices,
as in the proof
of Theorem \ref{lem22}, we find that $p_k^n({\bf x})$ is a factor of
\be
\label{eqn:bigdet}
{\rm det} \bigl(\, d \cdot I_{m^k}  \,\, + \,\,
{\bf U}_1({\bf x}) \oplus {\bf U}_2({\bf x}) \oplus
\cdots \oplus {\bf U}_k({\bf x})  \,\bigr).
\ee
However, there are many extraneous factors. They
are powers of the irreducible polynomials
that define the $k'$-ellipsoids whose foci
are subsets of $\{{\bf u}_1, {\bf u}_2, \ldots, {\bf u}_k\}$.

There is a standard choice for the matrices $\,{\bf U}_i({\bf x}) \,$ which
is symmetric with respect to permutations of the $n$ coordinates.
Namely, we can take $m = n+1$ and
$$ {\bf U}_i ( {\bf x} ) \quad = \quad
\bbm 0 & x_1 {-} u_{i1} & x_2 {-} u_{i2} & \cdots & x_n {-} u_{in} \\
\, x_1 - u_{i1} & 0                 &       0              &  \cdots &         0         \\
\, x_2 - u_{i2} & 0                 &       0              &  \cdots &         0         \\
\,\vdots    & \vdots         &    \vdots        &  \ddots  &  \vdots      \\
\,x_n - u_{in} & 0                 &       0              &  \cdots &         0         \\
 \ebm
$$ However, in view of the extraneous factors in (\ref{eqn:bigdet}),
it is desirable to replace these by matrices of smaller size $m$,
possibly at the expense of having additional nonzero eigenvalues.
It is not hard to see that $n$ is the smallest possible matrix size in a
symmetric determinant representation.
To see this, assume $d^2-\sum_{j=1}^n (x_j-u_{ij})^2 = \det(A_0+x_1A_1+\cdots+x_nA_n)$
where $A_i$ are all constant symmetric matrices of size $m\times m$.
Then it must hold $m\geq n$.
Otherwise, if $m<n$, the first row of $A_0+x_1A_1+\cdots+x_nA_n$
would vanish on some affine subspace of $\re^n$
and so did $d^2-\sum_{j=1}^n (x_j-u_{ij})^2$,
but this is not possible since $\{x\in \re^n:\, d^2-\sum_{j=1}^n (x_j-u_{ij})^2=0\}$
is compact.
However, if we allow extraneous factors,
then the smallest possible value of $m$ might drop.
These questions are closely related to find
the {\em determinant complexity} of a given polynomial,
as discussed in \cite{Mignon}. Note that in the applications
to complexity theory considered there, the matrices
of linear forms need not be symmetric.

\section{Open questions and further research}

The $k$-ellipse is an appealing example of an object from algebraic
geometry.  Its definition is elementary and intuitive, and yet it
serves well in illustrating the intriguing interplay between algebraic
concepts and convex optimization, in particular semidefinite
programming.  The developments presented in this paper motivate many
natural questions. For most of these, to the best of our knowledge, we
currently lack definite answers. Here is a short list of open problems
and possible topics of future research.

\paragraph{Singularities and genus}
Both the circle and the ellipse are rational curves, i.e., have genus
zero. What is the genus of the (projective) algebraic $k$-ellipse?
 The first values, from
$k=1$ to $k=4$, are 0,0,3,6. What is the formula for the genus in general?
The genus is related to the class of the curve, i.e.~the degree
of the dual curve, and this number is the algebraic degree \cite{NRS}
of  the problem (\ref{ellipseSDP}).
Moreover, is there a nice geometric characterization of all
(complex) singular points of the algebraic $k$-ellipse?

\paragraph{Algebraic degree of the Fermat-Weber point}
The Fermat-Weber point  $(x_*,y_*)$ is the unique solution of an
algebraic optimization problem, formulated in
(\ref{FWSDP}) or (\ref{betterFWSDP}),
 and hence it has a well-defined algebraic
degree over $\Q(u_1,v_1,\ldots,u_k,v_k)$. However, that
algebraic degree
will depend on the  combinatorial configuration of the
foci. For instance, in the case $k=4$ and foci forming a convex
quadrilateral, the Fermat-Weber point lies in the intersection of the
two diagonals \cite{Bajaj}, and therefore its algebraic degree is
equal to one. What are the possible values for this degree? Perhaps a
possible approach to this question would be to combine the results and
formulas in \cite{NRS} with the semidefinite characterizations
obtained in this paper.

\paragraph{Reduced SDP representations of rigidly convex curves}
A natural question motivated by our discussion in
Section~\ref{sec:semidefinite} is how to systematically produce
minimal determinantal representations for a rigidly convex curve, when
a non-minimal one is available. This is likely an easier task than
finding a representation directly from the defining polynomial, since
in this case we have a certificate of its rigid convexity.

Concretely, given real symmetric $n \times n$ matrices $A$ and $B$ such
that
\[
p(x,y) = \det (A \cdot x + B \cdot y + I_n)
\]
is a polynomial of degree $r < n$, we want to produce $r \times r$ matrices $\tilde A$
and $\tilde B$ such that
\[
p(x,y) = \det (\tilde A \cdot x + \tilde B \cdot y + I_r).
\]
The existence of such matrices is guaranteed by the results in
\cite{HV,LPR}. In fact, explicit formulas
in terms of theta functions of a Jacobian
variety are presented in \cite{HV}.
 But isn't there a simpler algebraic  construction in this special case?

\paragraph{Elimination in semidefinite programming}
The projection of an algebraic variety is (up to Zariski
closure, and over an algebraically closed field)
 again an algebraic variety. That projection can be computed using
elimination theory or Gr\"{o}bner bases.  The projection of a polyhedron
into a lower-dimensional subspace is a polyhedron. That projection can be
computed using Fourier-Motzkin elimination. In contrast
to these examples, the  class of feasible sets of semidefinite
programs is not closed under projections. As a simple concrete
example, consider the convex set
\[
\left\{ (x,y,t) \in \re^3 \; : \:
\bbm 1 & x-t \\ x-t & y \ebm \succeq 0, \quad t \geq 0 \right\}.
\]
Its projection onto the $(x,y)$-plane is a convex set that
is not rigidly convex, and hence cannot be expressed as
$\{(x,y) \in \re^2: A x + B y + C \succeq 0 \}$. In fact, that
projection is not even basic semialgebraic.
In some cases, however, this closure property nevertheless does hold.
We saw this for the projection
that transforms the representation
(\ref{betterFWSDP}) of the $k$-ellipse
to the representation (\ref{FWSDP}).
 Are there general conditions that ensure the semidefinite
representability of the projections? Are there
situations where the projection does not lead to an exponential
blowup in the size of the representation?

\paragraph{Hypersurfaces defined by eigenvalue sums}

Our construction of the (weighted) $k$-ellipsoid
as the determinant of a certain tensor sum has
the following natural generalization.
Let ${\bf U}_1({\bf x}) , {\bf U}_2({\bf x}), \ldots,
{\bf U}_k({\bf x})$ be any symmetric
$m \times m$-matrices whose entries
are affine-linear forms in ${\bf x} = (x_1,x_2,\ldots, x_n)$.
Then we consider the polynomial
\be
\label{specialhypersurface} p({\bf x}) \quad = \quad
{\rm det} \bigl(
{\bf U}_1({\bf x}) \oplus {\bf U}_2({\bf x}) \oplus \cdots
\oplus {\bf U}_k( {\bf x} ) \bigr) .
\ee
We also consider the corresponding rigidly convex set
$$
\bigl\{\, {\bf x} \in \re^n \,\,: \,\,
{\bf U}_1({\bf x}) \oplus {\bf U}_2({\bf x}) \oplus \cdots
\oplus {\bf U}_k( {\bf x} ) \, \succeq \, 0 \,
\bigr\}. $$
The boundary of this convex set is a hypersurface
whose Zariski closure is the set of zeroes of the
polynomial $p({\bf x})$.
It would be worthwhile to study the hypersurfaces
of the special form (\ref{specialhypersurface}) from
the point of view of computational algebraic geometry.

These hypersurfaces specified by eigenvalue sums
of symmetric matrices of linear forms have a natural generalization
in terms of {\em resultant sums} of hyperbolic polynomials.
For concreteness,
let us take $k=2$. If $p({\bf x})$ and
$q({\bf x})$ are hyperbolic polynomials in $n$ unknowns,
with respect to a common direction ${\bf e}$
in $\re^n$, then the polynomial
$$
 (p \oplus q)({\bf x}) \quad := \quad {\rm Res}_t
\bigl( p({\bf x} - t {\bf e}) , q({\bf x} + t {\bf e} ) \bigr) $$ is
also hyperbolic with respect to ${\bf e}$. This construction mirrors
the operation of taking Minkowski sums in the context of convex
polyhedra, and we believe that it is fundamental for future studies
of hyperbolicity in polynomial optimization \cite{LPR, Ren}.

\bigskip \bigskip
\bigskip \bigskip

\noindent {\bf Acknowledgements.}
We are grateful to the IMA in Minneapolis for hosting us during our
collaboration on this project.  Bernd Sturmfels was partially
supported by the U.S.~National Science Foundation (DMS-0456960).
Pablo A.~Parrilo was partially supported by AFOSR MURI subaward
2003-07688-1 and the Singapore-MIT Alliance.

\end{document}